\documentclass[conference]{IEEEtran}
\usepackage{fullpage}
\usepackage{amsmath}
\usepackage{amssymb}
\usepackage{marginnote}
\usepackage[final]{graphicx}
\usepackage{color}
\usepackage{mdwlist}
\usepackage{cite}
\newtheorem{theorem}{Theorem}
\newtheorem{lemma}[theorem]{Lemma}
\newtheorem{example}[theorem]{Example}
\newtheorem{rmk}[theorem]{Remark}
\newtheorem{defn}[theorem]{Definition}

\title{Link Biased Strategies in Network Formation Games}
\author{
\IEEEauthorblockN{Shaun Lichter and Terry Friesz}
\IEEEauthorblockA{
	Dept. Industrial and Manufacturing Engineering\\
	Penn State University\\
	E-mail: \{\texttt{tlf13}, \texttt{sml310}\}\texttt{@psu.edu}
} \and
\IEEEauthorblockN{Christopher Griffin}
\IEEEauthorblockA{
	Applied Research Laboratory\\
	Penn State University\\
	E-mail: \texttt{griffinch@ieee.org}
}
}
\date{March 15, 2011}
\begin{document}
\maketitle

\begin{abstract}
We show a simple method for constructing an infinite family of graph formation games with link bias so that the resulting games admits, as a \textit{pairwise stable} solution, a graph with an arbitrarily specified degree distribution. Pairwise stability is used as the equilibrium condition over the more commonly used Nash equilibrium to prevent the occurrence of ill-behaved equilibrium strategies that do not occur in ordinary play. We construct this family of games by solving an integer programming problem whose constraints enforce the terminal pairwise stability property we desire.
\end{abstract}

\section{Introduction}
The Network Science community has largely dedicated its efforts to the exposition and analysis of topological properties that occur in several real-world networks (e.g., scale-freeness \cite{barabasi1999a, newman2003, dorogovtsev2002, albert2002}). Recently, there has been interest in showing that these topological properties may arise as a result of optimization, rather than some immutable physical law \cite{doyle2005}. As Doyle et al. \cite{doyle2005} point out, various networks such as communications networks and the Internet are designed by engineers with some objectives and constraints. While it is true that there is often not a single designer in control of the entire network, the network does not \emph{naturally} evolve without the influence of \emph{designers}.  In each application, the network structure must be feasible with respect to some physical constraints corresponding to the tolerances and specifications of the equipment used in the network.  For example, in a system such as the world wide web, a single web-page might have billions of connections, however it is not possible for a node to have such a degree in many other applications, such as collaboration or road networks.  Certainly the structure of the network has a significant impact on its [the network's] ability to function, its evolution, and its robustness.  However network structures often arise as a result (locally) of optimized decision made by a single agent or multiple competitive or cooperative agents, who take network structure and function into account as a part of a collection of constraints and objectives.

Recently, network formation has been modeled from a game theoretic perspective \cite{myerson1977,jackson1996,dutta1997,goyal2003}, and in \cite{LichterGriffinFriesz2011}, it was shown that there exists games that result in the formation of a stable graph with an arbitrary degree sequence $\mathbf{k} = (k_1,\dots,k_n)$.  These models require that the players have similar objectives that are convex with minima near the desired $k_{i}$ values.  These assumptions were necessary to show that a game can be constructed that admits a stable graph with an arbitrary degree sequence.  However, this assumption is limiting in its modeling power because players (usually) do not have an exact number of links that they desire nor do they (usually) have an objective function precisely specifying this desire. Instead the $k_i$ arise endogenously as a result of other factors.  In this paper, we present a model incorporating a player's link bias - preference of one link over another.  The incorporation of link bias allows the game to  result in a stable graph of arbitrary degree without requiring the degree sequence to be precisely coded into the game.

\section{Model}
Let $N=\{1,2, \ldots n\}$ be a set of nodes. We will assume $n$ is fixed for the remainder of this paper. A link between two nodes is any subset of size 2 of set $N$. A graph $g$ is any set of links and the complete graph $g^c$ is the set of all size two subsets of $N$. The set $G$ is composed of all graphs over the node set $N$, that is, $G=\{g : g \subseteq g^{c} \}$.

In a network formation game, each node is a player. Link bias is introduced by assuming player $i$ has a cost function $f_{i}:G \rightarrow \mathbb{R}$.  For this paper, we assume a linear cost function:
\begin{equation}
f_{i}(g)=\sum_{j}c_{ij}x_{ij}
\end{equation}
where
\begin{equation}
x_{ij} =
\begin{cases}
1 & \text{there is a link between $i$ and $j$ in $g$}\\
0 & \text{else}
\end{cases}
\label{eqn:xdef}
\end{equation}
A player's strategy is to determine to which nodes to link in order minimize cost (or maximize payoff). Following \cite{jackson1996}, a link between players $i$ and $j$ exists if and only if the two players decide to link. That is, a player may \textit{unilaterally} reject a link. This is consistent with friending policies on Facebook or linking policies on LinkedIn.

\subsection{Stability}
The value of a graph $g$ is the total value produced by agents in the graph; we denote the value of a graph as the function $v:G \rightarrow \mathbb{R}$ and the set of of all such value functions as $V$. An allocation rule $Y:V \times G \rightarrow \mathbb{R}^{n}$ distributes the value $v(g)$ among the agents in $g$.  Denote the value allocated to agent $i$ as $Y_{i}(v,g)$.  Since, the allocation rule must distribute the value of the network to all players, it must be \textit{balanced}; i.e., $\sum_{i} Y_{i}(v,g)=v(g)$ for all $(v,g) \in V \times G$.  The allocation rule governs how the value is distributed and thus makes a significant contribution to the model. In the sequel we define $Y_{i}(v,g) = -f_i(g)$ and define the value function $v$ \textit{implicitly} to ensure that it is balanced. Jackson and Wolinksy use \textit{pairwise stability} to model stable networks without the use of Nash equilibria \cite{jackson1996}.

\begin{defn}
A network $g$ with value function $v$ and allocation rule $Y$ is pairwise \textit{stable} if (and only if):
\begin{enumerate*}
\item for all $ij \in g$, $Y_{i}(v,g) \geq Y_{i}(v,g - ij)$ and
\item for all $ij \not\in g$, if $Y_{i}(v,g+ij) > Y_{i}(v,g)$, then $Y_{j}(v,g+ij) < Y_{j}(v,g)$
\end{enumerate*}
\end{defn}
Define $Y_{i}(v,g)=-f_{i}(g)$.  In this model, player $i$ will benefit from linking with any player $j$ whenever $c_{ij}<0$ and player $i$ will be penalized for linking with any player $j$ whenever $c_{ij}>0$.  We may have an unspecified behavior when $c_{ij} = c_{ji} = 0$. In this case it will neither help nor hinder either player to establish a link. To remove this possibility, we may assume link parsimony. That is, we will assume that a link is established if and only if both players benefit in some way. The stability condition becomes:
\begin{defn}
A network $g$ with $Y_{i}(v,g)=-f_{i}(g)$ is pairwise \textit{stable} if (and only if):
\begin{enumerate*}
\item for all $ij \in g$, $c_{ij} < 0$ and $c_{ji} < 0$
\item for all $ij \not\in g$, if $c_{ij} \geq 0$ then $c_{ji} \leq 0$
\end{enumerate*}
\end{defn}

We can construct a cost matrix $\mathbf{C}$ via an optimization problem such that the resulting graph formation game has as a (pairwise) stable solution a graph with an arbitrary degree sequence $\mathbf{k} = (k_1,\dots,k_n)$.
\begin{defn}
Define $\boldsymbol{\psi}$ as a matrix of $0-1$ values indicating a player's interest in a particular link. That is:
\begin{equation*}
\boldsymbol{\psi}_{ij} = \begin{cases}
1 & \text{if player } i \text{ can benefit from link } ij\\
0 & \text{if player } i \text{ cannot benefit from link } ij
\end{cases}
\label{defn: psi}
\end{equation*}


Specifically, $\boldsymbol{\psi}$ is the boolean mapping of $\mathbf{C}$:

\begin{equation*}
\boldsymbol{\psi}_{ij} = \begin{cases}
1 & \text{if } c_{ij}<0\\
0 & \text{if } c_{ij} \geq 0
\end{cases}
\end{equation*}
\end{defn}
\begin{rmk} To each graph $g$, we may associate a vector of binary variables $\mathbf{x} = \langle{x_{ij}}\rangle$ where $x_{ij}$ is defined in Equation \ref{eqn:xdef}. For simplicity we will write $\mathbf{x}$ for the graph it represents when it is convenient.
\end{rmk}
\begin{lemma}
\label{lem: stability}
A graph  $\mathbf{x} = \langle{x_{ij}}\rangle$ is a symmetric stable graph if and only if it meets the following constraints:
\begin{equation}
\label{eqn: stability constraints0}
\begin{aligned}
x_{ij}&=x_{ji} \quad \forall \; ij\\
\psi_{ij}+\psi_{ji}-1 &\leq x_{ij} \quad \forall \; ij \\
x_{ij} &\leq \psi_{ij}  \quad \forall \; ij \\
x_{ij} &\leq \psi_{ji}  \quad \forall \; ij \\
x_{ij}, \psi_{ij} &\in \{0,1\} \quad \forall \; ij
\end{aligned}
\end{equation}
\end{lemma}
\begin{IEEEproof}
By definition, the first constraint ensures the graph is symmetric.  The second, third, and fourth constraints ensure stability.  If $\psi_{ij}=\psi_{ji}=1$, then by the definition of $\boldsymbol{\psi}$, both $c_{ij}<0$ and $c_{ji}<0$ and the second constraint forces $x_{ij}=1$ as stability requires.  Alternatively, if either $\psi_{ij}=0$ or $\psi_{ji}=0$ or both, then by the definition of $\boldsymbol{\psi}$, $c_{ij} \geq 0$ or $c_{ji} \geq 0$ or both and the third and fourth constraint force $x_{ij}=0$ as stability requires.  Lastly, if $x_{ij}=0$, then the second constraint implies that either $\psi_{ij}=0$ ($c_{ij} \geq 0$) and player $i$ will veto link $ij$ or $\psi_{ji}=0$ ($c_{ji} \geq 0$) and $j$ will veto link $ij$, as stability requires or $\psi_{ij} = \psi_{ji} = 0$ and both $c_{ij} \geq 0$ and $c_{ji} \geq 0$.

Now, suppose there is a symmetric stable graph $\mathbf{x}$ that violates one of these constraints.  If $\mathbf{x}$ violates the first constraint, then it is not a symmetric graph.  If the graph $\mathbf{x}$ violates the second constraint, then there exists a link $ij$ such that $\psi_{ij}=1$, $\psi_{ji}=1$, but $x_{ij}=0$.  This implies that $c_{ij}<0$ and $c_{ji}<0$ which implies that both player $i$ and $j$ can benefit from link $ij$ and this requires by stability that $x_{ij}=1$.  However, this is a contradiction.  If $\mathbf{x}$ violates the third or fourth constraint, then there exists a link $ij$ such that $x_{ij}=1$, but either $\psi_{ij}=0$ ($c_{ij} \geq 0$) or $\psi_{ji}=0$ ($c_{ji} \geq 0$), but this would violate stability and hence be a contradiction.  There cannot be a symmetric pairwise stable graph $\mathbf{x}$ (and thus is a solution to the graph formation game with link bias matrix $\mathbf{C}$) that violates one of these constraints.
\end{IEEEproof}

If these constraints are consistent (i.e., the feasible region is non-empty), then any feasible $\boldsymbol{\psi}$ for the constraints above can be used to generate a cost matrix $\mathbf{C}$. In fact, it can be used to generate an infinite number of cost matrices, the simplest one given by:
\begin{equation}
\mathbf{C}_{ij} = \begin{cases}
-1 & \text{if }\psi_{ij} = 1\\
1 & \text{if }\psi_{ij} = 0
\end{cases}
\label{defn:C}
\end{equation}


\begin{lemma}
\label{lem: stabilityK}
Suppose $\mathbf{C}$ is a cost matrix and $f_1,\dots,f_N : G \rightarrow \mathbb{R}$ are player cost functions with $Y_i(v,\mathbf{x}) = -f_i(\mathbf{x}) = -\sum_{j}c_{ij}x_{ij}$, which define a graph formation game (for a graph $\mathbf{x}$). Then the graph $\mathbf{x} = \langle{x_{ij}}\rangle$ is symmetric, pairwise stable (for the game) and has degree sequence $\mathbf{k}=(k_{1},k_{2},\ldots, k_{n})$ if and only if it meets the following constraints:
\begin{equation}
\begin{aligned}
\sum_{j \neq i} x_{ij}&=k_{i} &\quad\forall \; i  \\
x_{ij}&=x_{ji} &\quad \forall \; ij\\
\psi_{ij}+\psi_{ji}-1 &\leq x_{ij} &\quad \forall \; ij \\
x_{ij} &\leq \psi_{ij}  \quad \forall \; ij \\
x_{ij} &\leq \psi_{ji}  \quad \forall \; ij \\
\psi_{ij} &\in \{0,1\} &\quad \forall \; ij
\end{aligned}
\label{eqn: stability constraints1}
\end{equation}
\end{lemma}
\begin{IEEEproof}
The first constraint ensures that the graph has degree sequence $\mathbf{k}=(k_{1},k_{2},\ldots, k_{n})$. The remainder of the proof follows from Lemma \ref{lem: stability}.
\end{IEEEproof}

\subsection{Construction of Cost Matrix via Optimization}
In the event that there is no feasible $\boldsymbol{\psi}$, we may solve an optimization problem to find close solutions in the $\ell_1$ metric on vector $\mathbf{x}$ by pricing out the first constraints in Problem \ref{eqn: stability constraints1}:
\begin{equation}
\label{model: optimization 1A}
\begin{aligned}
\min\;\; & \sum_{i} \left| \sum_{j \neq i} x_{ij}-k_{i} \right|\\
s.t.\; & x_{ij} - x_{ji} = 0 \quad \forall i < j\\
& \psi_{ij}+\psi_{ji}-1 \leq x_{ij} \quad \forall \; ij\\
x_{ij} &\leq \psi_{ij}  \quad \forall \; ij \\
x_{ij} &\leq \psi_{ji}  \quad \forall \; ij \\
& x_{ij},\psi_{ij} \in \{0,1\} \quad \forall i,j
\end{aligned}
\end{equation}
We may reformulate this optimization problem as:
\begin{equation}
\label{model: optimization 1B}
\begin{aligned}
\min\;\; & \sum_{i} e_{i}\\
s.t.\; &  \left| \sum_{j \neq i} x_{ij}-k_{i} \right| \leq e_{i}\\
&x_{ij} - x_{ji} = 0 \quad \forall i < j\\
& \psi_{ij}+\psi_{ji}-1 \leq x_{ij} \quad \forall \; ij\\
x_{ij} &\leq \psi_{ij}  \quad \forall \; ij \\
x_{ij} &\leq \psi_{ji}  \quad \forall \; ij \\
& x_{ij},\psi_{ij} \in \{0,1\} \quad \forall i,j\\
\end{aligned}
\end{equation}
The non-linear constraints can be transformed into equivalent linear constraints with the result being an integer linear programming problem.
\begin{equation}
\label{model: optimization 1C}
\begin{aligned}
\min\;\;  &\sum_{i} e_{i} &\\
s.t.\;   &\sum_{j \neq i} x_{ij}-k_{i}  \leq e_{i}\\
 - &\sum_{j \neq i} x_{ij}+k_{i}  \leq e_{i}\\
 & \psi_{ij}+\psi_{ji}-1 \leq x_{ij} \quad \forall \; ij\\
 x_{ij} &\leq \psi_{ij}  \quad \forall \; ij \\
x_{ij} &\leq \psi_{ji}  \quad \forall \; ij \\
&x_{ij} - x_{ji}= 0 \quad \forall i < j\\
& x_{ij},\psi_{ij} \in \{0,1\} \quad \forall i,j\\
\end{aligned}
\end{equation}
Once a solution $\boldsymbol{\psi}$ is found to this optimization problem, we can construct a cost matrix $\mathbf{C}$ using (e.g.,) Expression \ref{defn:C}.  Any matrix $\mathbf{C}$ constructed from such a solution $\boldsymbol{\psi}$ will result in a graph as close as possible to degree sequence $\mathbf{k} = (k_1,\dots,k_n)$ in the $\ell_1$ metric. The nature of Problem \ref{model: optimization 1C} ensures the following theorem:

\begin{theorem} \label{thm: general arbitrary dist} For a given degree sequence $\mathbf{k}$ any cost matrix $\mathbf{C}$ derived as a solution to Problem \ref{model: optimization 1C} and by (e.g.,) Expression \ref{defn:C} will yield a graph formation game that itself admits, as a pairwise stable solution, a graph $g$ whose degree sequence is as close as possible to $\mathbf{k}$ in the $\ell_1$ metric.
\end{theorem}
\begin{IEEEproof}
As shown in Lemma \ref{lem: stability}, any graph that satisfies the third, fourth, and fifth constraints is stable. Consider the simplified problem:
\begin{equation}
\label{model: manhattan min}
\begin{aligned}
\min\;\;  & e &\\
s.t.\;   &x-k  \leq e \\
 - &x+k  \leq e\\
\end{aligned}
\end{equation}
It is easy to see that that the solution to Problem (\ref{model: manhattan min}) yields the closest $x$ to $k$ in the $\ell_1$ metric.  From here the proof is straight forward.

Denote the solution to Problem (\ref{model: optimization 1C}) as $\mathbf{x}^{\ast}$ and the objective function value at this solution as $\sum_{i} e^{\ast}_{i}$.  Suppose there is a graph $\mathbf{x}$ that is both stable and has a degree sequence closer to $\mathbf{k} = (k_1,\dots,k_n)$ than $\mathbf{x}^{\ast}$ in the $\ell_1$ metric.  This implies that:
\begin{equation*}
\sum_{i} \sum_{j \neq i} |x_{ij}-k_{i}| < \sum_{i} \sum_{j \neq i} |x^{\ast}_{ij}-k_{i}| \leq \sum_{i} e^{\ast}_{i}
\end{equation*}
and $\mathbf{x}$ satisfies the other constraints for stability.  This implies that $\mathbf{x}^{\ast}$ is not the optimal solution to (\ref{model: optimization 1C}), which is a contradiction.
\end{IEEEproof}

\section{Numerical Example}
\begin{example}
Suppose that we want the degree sequence of a stable graph that results from playing the game described in Theorem \ref{thm: general arbitrary dist} to have a power law degree distribution with the degree of each player ($k_i$) given in the following table:
\begin{center}
\begin{tabular}{|c|c|}
\hline
Node(s) & $k_{i}$ \\
\hline
1-29 & 1\\
\hline
30-33 & 2\\
\hline
34 & 3\\
\hline
35 & 4\\
\hline
\end{tabular}
\end{center}
The distribution of values of $k_{i}$ forms an approximate (with rounding to integers) power law distribution.  This is illustrated in Figure \ref{fig:subfig2}:
\begin{figure}[htbp]
\centering
\includegraphics[scale=0.4]{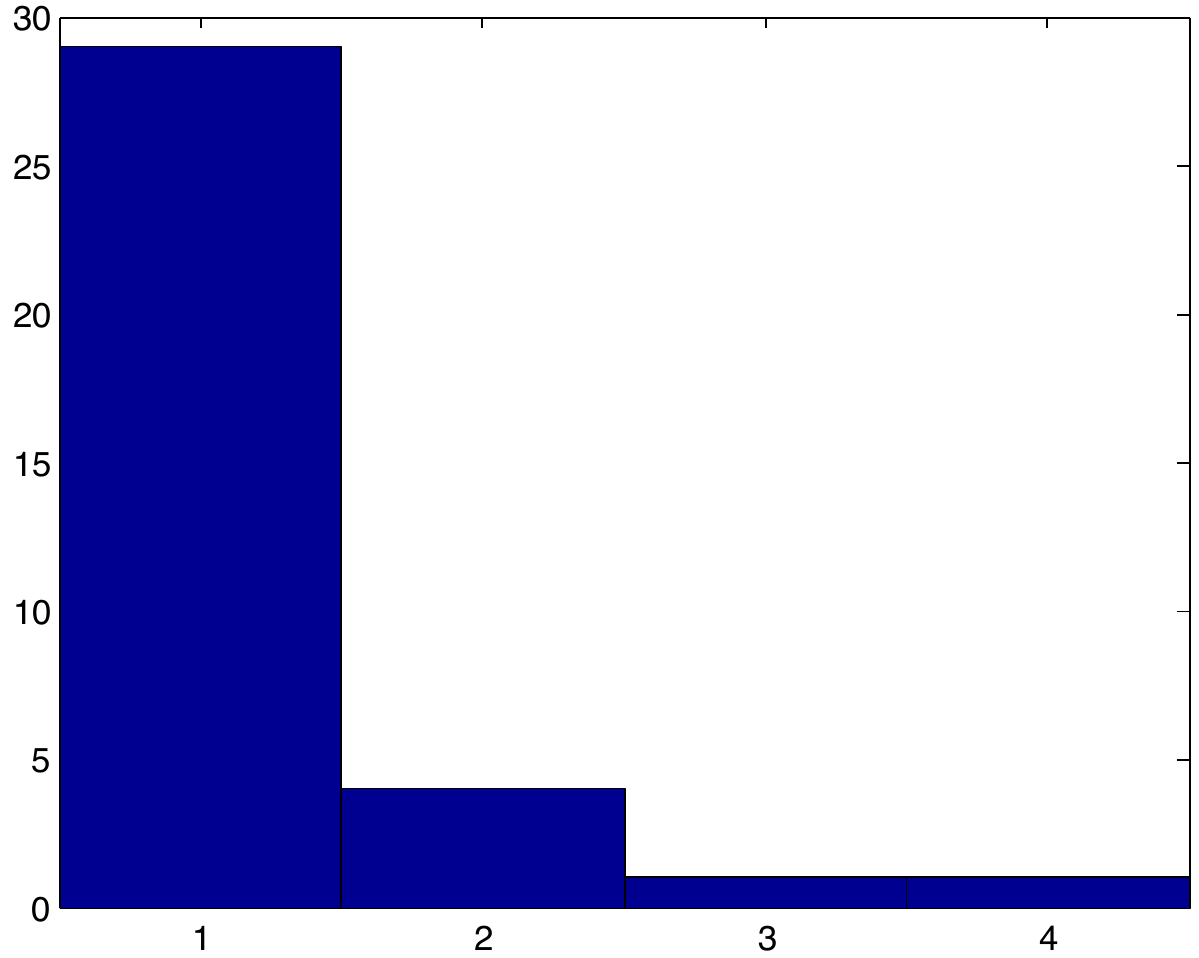}
\caption{The empirical desired distribution of the degrees of the players. This degree distribution follows an approximate power law distribution.}
\label{fig:subfig2}
\end{figure}
We can now solve the optimization Problem \ref{model: optimization 1C} to find the stable graph with degree sequence as close (in Manhattan Norm) as possible to degree sequence $\mathbf{k}$.  The resulting graphic solution is shown in Figure \ref{fig:X}.
\begin{figure}[ht]
\centering
\includegraphics[scale=0.60]{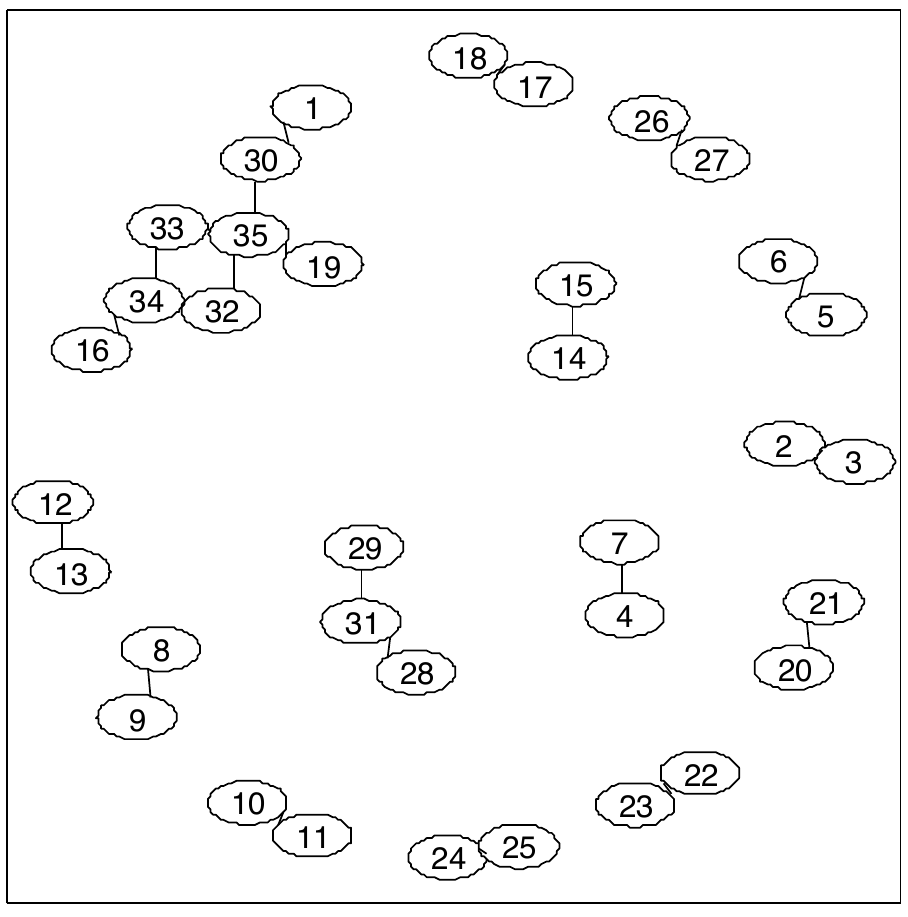}
\caption{Collaboration Network for 35 nodes}
\label{fig:X}
\end{figure}
The degree distribution of the graphic solution is compared to the distribution of the $k_{i}$ values in Figure \ref{fig:X2}.
\begin{figure}[ht]
\centering
\includegraphics[scale=0.5]{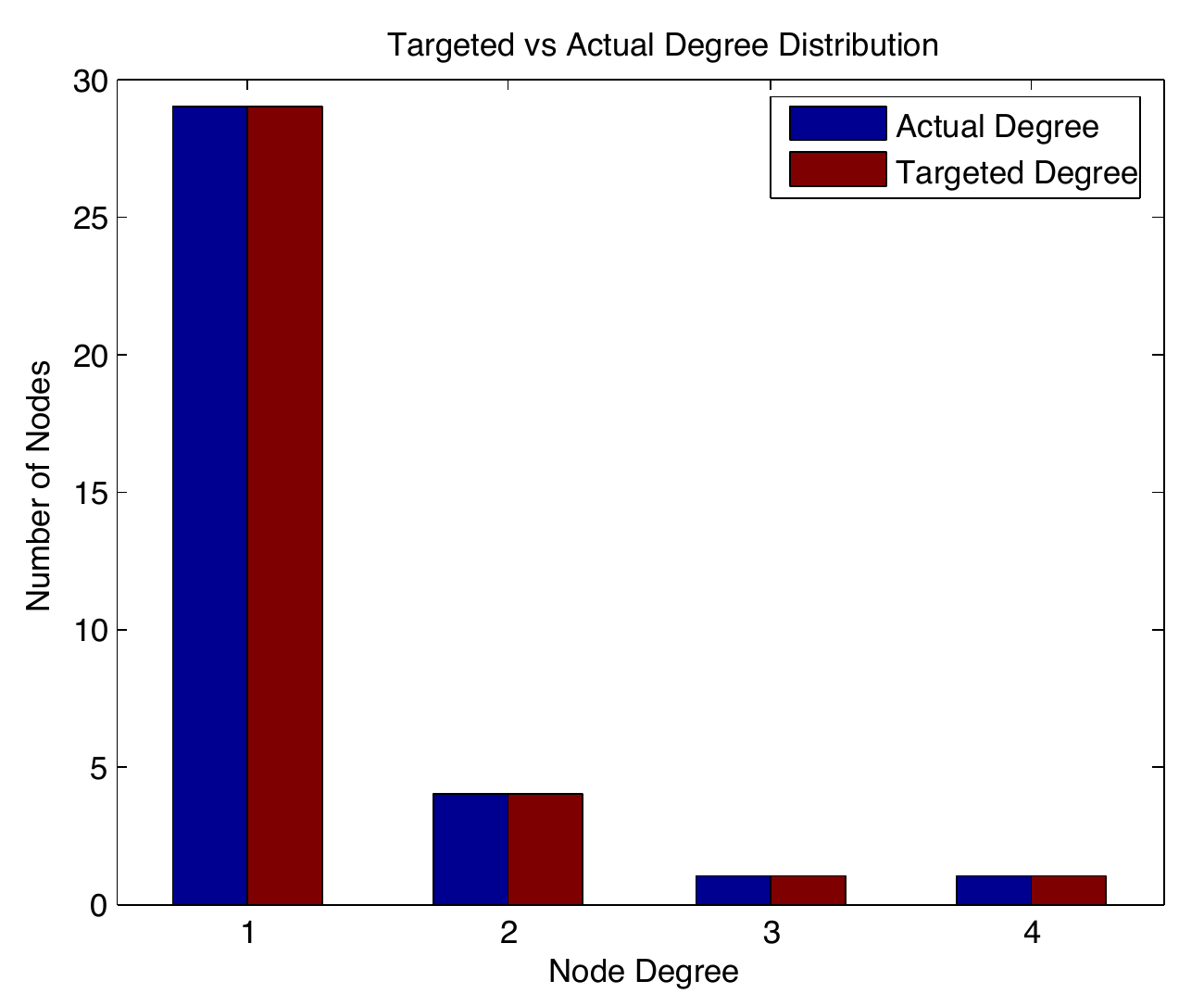}
\caption{Degree Distribution Comparison}
\label{fig:X2}
\end{figure}

\begin{table}[htbp]
\centering
\begin{tabular}{|cc|cc|cc|cc|}
  \hline
  i & j & i & j & i & j & i & j\\
  \hline
  1 &  30 &  2 &  3 &  3 &  2 & 34 & 33\\
  10 & 11 & 20 & 21 & 30 &  1 & 35 & 19\\
  11 & 10 & 21 & 20 & 30 & 35 & 35 & 30\\
  12 & 13 & 22 & 23 & 31 & 28 & 35 & 32\\
  13 & 12 & 23 & 22 & 31 & 29 & 35 & 33\\
  14 & 15 & 24 & 25 & 32 & 34 & 4 & 7\\
  15 & 14 & 25 & 24 & 32 & 35 & 5 & 6\\
  16 & 34 & 26 & 27 & 33 & 34 & 6 & 5\\
  17 & 18 & 27 & 26 & 33 & 35 & 7 & 4\\
  18 & 17 & 28 & 31 & 34 & 16 & 8 & 9\\
  19 & 35 & 29 & 31 & 34 & 32 & 9 & 8\\
  \hline
\end{tabular}
\caption{Table illustrating elements of $\mathbf{C}$ that have value $1$. All other entries are $-1$.}
\label{tab:Entries}
\end{table}
Now, we can construct $\mathbf{C}$ from $\boldsymbol{\psi}$. If $\psi_{ij}=1$ (the pair is listed in Table \ref{tab:Entries}), then $c_{ij}=-1$, otherwise $c_{ij}=1$.
\end{example}

\section{Linking Portfolio with Resource Constraint}
The model presented in \cite{LichterGriffinFriesz2011} assumes a player has a specific desired degree.  Alternatively, the model presented in the previous section presumes instead that a player links to players that can benefit them, but the number of such links is unlimited. In reality players link to other players who benefit them, but they are also constrained by time and other resources that they must consider, which prohibits them from befriending everyone they could benefit from.  Consider player $i$ may solve an optimization problem with form:
\begin{equation}
\label{model: portfolio optimization}
\begin{aligned}
\max \; \; & f_i(\mathbf{x})\\
s.t.\; & \sum_{j} a_{ij}x_{ij} \leq b_{i}\\
& x_{ij} \in \{0,1\} \quad \forall i,j\\
\end{aligned}
\end{equation}
Each set of possible links $\mathbf{x}_{i}=(x_{i1},x_{i2},\ldots x_{in})$ is mapped to a ranking (which may have ties) by this math program.  Since a player now considers their entire strategy when making a linking decision, we cannot attribute the existence of a link or lack there of, to a specific characteristic of only that link ($c_{ij}$ or $a_{ij}$), rather it is a result of the entire vector of parameters.

Nonetheless, given a degree sequence $\mathbf{k}$, we may construct a resource matrix $\mathbf{A}$ and link cost matrix $\mathbf{C}$, such that the resulting graph is stable and has degree sequence $\mathbf{k}$.  It is rather straight forward to show this is possible using to Theorem \ref{thm: general arbitrary dist}.  We may consider $f_i(\mathbf{x})=-\sum_{j} c_{ij}x_{ij}$, $\mathbf{A}=0$, and large $\mathbf{b}$, which essentially makes the resource constraint void.  This model is then equivalent to the model in Theorem \ref{thm: general arbitrary dist}.

However, this model allows players to consider not only their limited resources, but the fact that some links require greater resources than others and some links provide greater benefits than others and these benefits and resource consumption may be asymmetric between two players for a given link.  Further, a general $f$ will allow a player to consider one player's link with greater weight than other or in fact a player may consider nonlinear benefit functions.  For example, a player may benefit from one subset of players more than another or there may be nonlinear effects between links with other players.  Nonetheless, for an arbitrary degree sequence $\mathbf{k}$, there is always a game that can be constructed with this model that has a stable graph with that degree sequence. A theoretical characterization of these possibilities is left for future work.
\newpage
\section{Conclusion}
With modest modeling characteristics, we may construct a network formation game that admits as a pairwise stable solution a graph with an arbitrary degree sequence. In fact, we may construct an entire family of games admitting such graphs as pairwise stable solutions. We can likely  construct games with other graph characteristics as well. Network Formation Games are capable of modeling a wide variety of graph structures.  Further, there is opportunity to add richness to the model with aspects such as the resource budget link portfolio model.  




\bibliographystyle{IEEEtran}
\bibliography{Biblio-Database2}

\end{document}